\documentclass[12pt,a4paper]{article}
\usepackage{a4wide,amsmath,amsthm,epsfig,graphicx,psfrag}
\usepackage{amsmath,amsthm}
\usepackage{amsfonts}
\usepackage{amssymb}
\usepackage{mathscinet}
\usepackage{enumerate,fontenc}
\usepackage[latin1]{inputenc}  

\newcommand{\url}[1]{}

\newcommand{\R}{{\mathbb R}}
\newcommand{\Ex}{\Bbb{E}}

\newcommand{\Hs}{\mathcal{H}_{\text{\bf vol}}}
\newcommand{\Hsw}{\mathcal{H}_{\text{\bf vol}}'}

\newcommand{\tsigma}{\tilde{\sigma}}

\newcommand{\ysd}{y^*}
\newcommand{\xss}{x^*_\sigma}

\newcommand{\yse}{\ysd_\eta}

\newcommand{\utsigma}{\utilde{\sigma}}
\newcommand{\uteta}{\utilde{\eta}}

\newtheorem{theorem}{Theorem}[section]
\newtheorem{proposition}[theorem]{Proposition}

\newtheorem{example}[theorem]{Example}
\newtheorem{remark}[theorem]{Remark}

\def\title#1{\hfil\break\hfil\break
\hfil\break\par\addvspace\baselineskip\noindent
\ignorespaces{\LARGE\bf#1}\hfil\break}

\def\author#1{\par\addvspace\baselineskip\noindent
\ignorespaces{\large\bf#1}}

\def\institute#1{\par\addvspace\baselineskip\noindent
\ignorespaces{\small#1}\hfil\break}

\begin{document}
\title{{\bf General Duality for Perpetual American Options }}
\author{Aur\'{e}lien Alfonsi$^{1,2}$ and Benjamin Jourdain$^1$}
\institute{$^1$CERMICS, project-team Mathfi, École des
     Ponts, ParisTech, 6-8
avenue Blaise Pascal, Cit\'{e} Descartes, Champs sur Marne,
77455 Marne-la-vall\'{e}e, France.\\
$^2$Institut f\"{u}r Mathematik, MA 7-4, TU Berlin, Stra\ss e des 17. Juni 136, 10623 Berlin, Germany.\\
e-mail : {\tt \{alfonsi,jourdain\}@cermics.enpc.fr } \\
\today }
\begin{abstract}
   In this paper, we investigate the generalization of the Call-Put
   duality equality obtained in \cite{AJ} for perpetual American options
   when the Call-Put payoff $(y-x)^+$ is replaced by $\phi(x,y)$. It
   turns out that the duality still holds under monotonicity and
   concavity assumptions on $\phi$. The specific analytical form of the
   Call-Put payoff only makes calculations easier but is not crucial
   unlike in the derivation of the Call-Put duality equality for European options.
   Last, we give some examples for which the optimal strategy is known explicitly. 
\end{abstract}

\section*{Introduction}

In~\cite{AJ}, we have obtained a Call-Put duality equality for perpetual American
options. More precisely, for an interest rate $r>0$, a dividend rate
$\delta\geq 0$ and a time-homogeneous local volatility function
$$(\Hs)\;\;\;\;\sigma :\R_+^*\rightarrow \R_+^*\mbox{ continuous and such that }\exists
\underline{\sigma},\overline{\sigma}\in \mathbb{R}^*_+,\;\forall x>0,\;\underline{\sigma} \le \sigma(x) \le
  \overline{\sigma},$$
let $(S^{x}_t,t \ge
0)$ denote the unique weak solution of the stochastic differential equation
\begin{equation}\label{EDS}
  dS^x_t=S^x_t((r-\delta)dt+\sigma(S^x_t)dW_t),\;S^x_0=x
\end{equation}
where $(W_t,t\ge 0)$ is a
standard Brownian motion.
We have proved the existence of another volatility function $\eta$
satisfying $(\Hs)$ such that \[\forall x,y>0,\;\underset{\tau \in
  \mathcal{T}_{0,\infty}} {\sup} \Ex \left[e^{-r\tau} (y-S^{x}_{\tau})^+
\right] = \underset{\tau \in
  \mathcal{T}_{0,\infty}} {\sup} \Ex \left[e^{- \delta \tau} (\overline{S}^{y}_{\tau}-x)^+
\right],\]
where $(\overline{S}^{y}_t,t \ge
0)$ denotes the weak solution to
\begin{equation}\label{EDS2}
 d\overline{S}^y_t=\overline{S}^y_t((\delta-r)dt+\eta(\overline{S}^y_t)dW_t),\;\overline{S}^y_0=y.
\end{equation}
Here, $\mathcal{T}_{0,\infty}$ denotes the set of stopping times with respect to the
usual natural filtration of the underlying.\\
Our primal goal was to generalize to American derivatives the Call-Put
duality equality 
\begin{equation}
  \mbox{ for }\eta\equiv\sigma,\;\forall T,x,y>0,\;\Ex \left[e^{-rT} (y-S^{x}_{T})^+
\right] = \Ex \left[e^{- \delta T} (\overline{S}^{y}_{T}-x)^+
\right]\label{eurdual}
\end{equation}\pagestyle{myheadings} \markboth{}{{\footnotesize General duality for Perpetual American Options }}
which holds in the European case. In
the perpetual American case, unless $\sigma$ is a constant (usual
Black-Scholes model), then $\eta$ is
different from $\sigma$. Then it is natural to wonder whether the
European and the perpetual American Call-Put dualities are similar in
nature. The European equality is equivalent to Dupire's formula~\cite{Dup} and,
to our knowledge, the existing proofs of both results rely heavily on the
following specifity of the payoff function $(y-x)^+$ : for 
fixed $y$ (resp. fixed $x$), in the distribution sense,
$\partial^2_x(y-x)^+=\delta_y(x)$
(resp. $\partial^2_y(x-y)^+=\delta_x(y)$). This property implies for instance that the second
order derivative of the European Put price with respect to the strike variable
$y$ is the actualized density of the underlying asset $S^x_T$. It
is not obvious at all that the equality \eqref{eurdual} could be generalized by replacing
respectively $(y-S^x_T)^+$ and $(\overline{S}^{y}_{T}-x)^+$ by $\phi(S^x_T,y)$ and $\phi(x,\overline{S}^{y}_{T})$.\\
In contrast, we are going to show that it is possible to generalize
the perpetual American duality to payoff functions $\phi(x,y)$ which
only share global properties with $(x-y)^+$. The loss of the specific
analytical expression only makes calculations more
complicated. 

More precisely, from now on, we assume that
$\phi:\R_+^*\times\R_+^*\rightarrow \R_+$ is a continuous function such
$\Phi=\{(x,y):\phi(x,y)>0\}\neq\emptyset$, $\phi$ is $C^2$ on $\Phi$ and such that
\begin{equation}
   \forall x,y\in \Phi,\;\partial_x\phi(x,y)<0,\;\partial_y\phi(x,y)>0,\;\partial^2_x\phi(x,y)\leq
  0\mbox{ and  }\partial^2_y\phi(x,y)\leq 0\label{hypphi}.
\end{equation}
Of course the function $(y-x)^+$ satisfies these assumptions. More
general examples are given in Section \ref{secduality}.\\
For $y>0$, let us define $X(y)=\inf \{ x>0, \ \phi(x,y)=0 \}$ with the convention
$\inf \emptyset = + \infty$. Thanks to~\eqref{hypphi},
we have $\{ x>0, \ \phi(x,y)=0 \}=\{x>0, x \ge X(y)\} $ and $0\le
X(y)<\infty$. Moreover, the
function $y \mapsto X(y)$ is nondecreasing. Let us also define $Y(x)=\inf\{y>0,
\phi(x,y)>0\}=\inf\{y>0, x <X(y)\}$.
As the pseudo-inverse of the nondecreasing function $X$, the function
$Y$ is nondecreasing. Finally,
$$\Phi=\{(x,y), \  \phi(x,y)>0 \}= \{(x,y), \ x < X(y) \}= \{(x,y), \ y > Y(x) \}.$$

We also make the following assumption weaker than $(\Hs)$ on
the volatility functions :
$$(\Hsw)\;\;\;\;\sigma :\R_+^*\rightarrow \R_+^*\mbox{ continuous and such that }\exists
\overline{\sigma}<+\infty,\;\forall x>0,\;\sigma(x)<\overline{\sigma}.$$
When $\sigma$ and $\eta$ satisfy $(\Hsw)$, then weak existence and
uniqueness hold for \eqref{EDS} and \eqref{EDS2} (see for example
Theorem 5.15 in \cite{KS}, using a log transformation). Let
\[ P_\sigma (x,y)=\underset{\tau \in
  \mathcal{T}_{0,\infty}} {\sup} \Ex \left[e^{-r\tau} \phi(S^{x}_{\tau},y)
\right] \text{ and } c_\eta (y,x)= \underset{\tau \in
   \mathcal{T}_{0,\infty}} {\sup} \Ex \left[e^{- \delta \tau}\phi(x,\overline{S}^{y}_{\tau})
 \right]. \]
where the notations $P$ and $c$ standing respectively for ``Put'' and ``Call'' are slightly abusive.


The paper is structured as follows. The first section is devoted to the pricing of
perpetual American options with payoff $\phi$. It turns out that, as in
the Call-Put case, for fixed strike $y>0$ (resp. $x>0$), there is a unique $x^*(y)$ (resp. $y^*(x)$) such that
$$\{x:P_\sigma(x,y)>\phi(x,y)\}=(x^*(y),+\infty)\mbox{
  resp. }\{y:c_\eta(y,x)>\phi(x,y)\}=(0,y^*(x)).$$
These exercise boundaries $x^*(y)$ and $y^*(x)$ are characterized by
some implicit equations involving $\phi$, and we prove that they solve explicit ODEs.
The second section deals with the duality result. We state a general result and, for
two specific families of payoff functions, we are able to find an explicit relation between dual
volatilities, as in the call-put case. 

\section{Pricing of the perpetual American options}

\subsection{Pricing formulas and exercise boundaries}

In this section we will use the approach of Beibel and Lerche \cite{BL}
to explicit the pricing functions $P_\sigma$ and $c_\eta$. As in \cite{AJ}, we will denote
by $f$ (resp. $g$) the unique, up to a multiplicative constant, positive
nonincreasing (resp. nondecreasing) solution of
\begin{eqnarray}
  \frac{1}{2} \sigma^2(x) x^2 f''(x)+ (r-\delta) x  f'(x) -  rf(x) &=&0, \ x>0 \label{EDf}\\
  \text{(resp. }  \frac{1}{2} \eta^2(x) x^2 g''(x)+ (\delta-r) x g'(x) -
  \delta g(x) &=&0, \ x>0 \text{ ).}\label{EDg}
\end{eqnarray}
Let us also recall that 
\begin{equation}\label{convex_fg}
  \forall x>0,\;f''(x)>0 \text{ and } g''(x)>0.
\end{equation}
This has been checked for example in~\cite{AJ} (Lemma 3.1) where
$\sigma$ is assumed to satisfy $(\Hs)$, but
the boundedness from below is not used in the proof.

\begin{proposition}\label{opt_put}Let us fix a strike $y>0$. If $X(y)=0$, then $\forall x>0,  P_\sigma(x,y)=0$. Otherwise there is a unique $\xss(y) \in
  (0,X(y))$ such that $\tau^P_x=\inf \{t \ge 0, S^x_t \le \xss(y) \}$ (convention $\inf
 \emptyset=+\infty$) is an optimal stopping time for $P_\sigma$ and:
 \begin{equation}\label{DEFx*}
  \forall x \le \xss(y), P_\sigma(x,y)=\phi(x,y), \ \forall x >
  \xss(y), P_\sigma(x,y)= \frac{\phi(\xss(y),y) }{
    f (\xss (y))}f (x) > \phi(x,y).  \hspace{0.1cm}
\end{equation}
In addition, we have 
\begin{equation}
   \frac{\phi (\xss(y),y)}{\partial_x \phi (\xss(y),y)}=
\frac{f(\xss(y))}{f'(\xss(y))}\label{defx*}
\end{equation} which implies the smooth-fit principle.
Last, the function $y\in \{ z: X(z)>0 \} \mapsto \xss(y)$ is~$\mathcal{C}^1$ and
satisfies the following ODE:
\begin{eqnarray}\label{EDOx*}
  \xss(y)'&=& \left[\partial^2_{xy}\phi (\xss(y),y) - \frac{\partial_x \phi
      (\xss(y),y)\partial_y \phi (\xss(y),y)}{\phi (\xss(y),y)}\right] \\
  && \times \frac{\xss(y)^2 \sigma^2(\xss(y))}{2[r\phi (\xss(y),y) +(\delta-r)\xss(y)
    \partial_x \phi (\xss(y),y) ]-\xss(y)^2 \sigma^2(\xss(y))\partial^2_x \phi
  (\xss(y),y)}.\nonumber \end{eqnarray}
It is strictly
increasing if one assumes moreover $\phi \partial ^2_{xy} \phi >\partial_x \phi\partial_y \phi $ on $\Phi$.
\end{proposition}
\begin{proof}
  For $x>0$, let $h(x)=\frac{\phi(x,y)}{f(x)}$. The function $h$ is
  nonnegative and we have
  $h(0^+)=0$ because $\underset{x\rightarrow 0^+}{\lim} f(x)= +\infty$ (see
  \cite{AJ}) and $\underset{x\rightarrow 0^+}{\lim}\phi(x,y)$ exists thanks to
  the monotonicity assumption and is finite thanks to the concavity
  assumption made in \eqref{hypphi}. We have also
  $h(x)=0$ for $x \ge X(y)$. Therefore the function $h$ reaches its
  maximum at some
  $\xss(y) \in (0,X(y))$. In particular we have $h'(\xss(y))=0$ which
  also writes $F(\xss(y),y)=0$ where the function
  $F(x,y)=\frac{\phi(x,y)}{\partial_x\phi(x,y)}-\frac{f(x)}{f'(x)}$ is
  defined on $\Phi$. This proves~\eqref{defx*}.\\Now
  since $\partial^2_x\phi(x,y)\le 0$ on
$\{x < X(y)\}$ and $\frac{f(x)f''(x)}{f'(x)^2}$ is a positive function (see
\eqref{convex_fg} for the convexity of $f$),
  $\partial_x F(x,y)= -\frac{\phi(x,y)\partial^2_x\phi(x,y)}{(\partial_x\phi(x,y))^2}
+\frac{f(x)f''(x)}{f'(x)^2}$ is positive on $(0,X(y))$ which ensures
uniqueness of $\xss(y)$.  The implicit function theorem, yields that $\xss$ is $\mathcal{C}^1$ in
the neighborhood of $y$ and $\xss(y)'=-\frac{\partial_y F(\xss(y),y)}{\partial_xF(\xss(y),y)}$. Since $\partial_y F(\xss(y),y)=\frac{\partial_x\phi(\xss(y),y)\partial_y\phi(\xss(y),y) -\phi(\xss(y),y)
  \partial^2_{xy}\phi(\xss(y),y)}{(\partial_x\phi(\xss(y),y))^2}$,
$\xss(y)'$ is positive if $\phi \partial^2_{xy} \phi >\partial_x \phi\partial_y \phi
$ on $\Phi$. From~\eqref{defx*} and the ODE~\eqref{EDf} satisfied by $f$, one gets
\begin{eqnarray*}
  \frac{f''(\xss(y))}{f'(\xss(y))}&=&\frac{2}{\xss(y)^2\sigma^2(\xss(y))}\left[r\frac{\phi (\xss(y),y)}{\partial_x \phi
      (\xss(y),y)}  +(\delta-r)\xss(y)\right]\end{eqnarray*}
so that we can express $\partial_xF(\xss(y),y)$ only with the derivatives of~$\phi$
and deduce~\eqref{EDOx*}.

When $x\ge \xss(y)$, the optimality of $\tau_x^P$ follows from the
arguments given in the proof of Theorem 1.4 \cite{AJ}. Let us now assume
that $x\in(0,\xss(y))$ and set $\tau^x_z=\inf\{t\geq 0:S^x_t=z\}$ for
$z>0$. Using the strong Markov property and the optimality result
when the initial spot is $\xss(y)$, then Fatou Lemma, we get for $\tau \in
  \mathcal{T}_{0,\infty}$, 
\[\Ex[e^{-r \tau}
  \phi(S^{x}_{\tau},y)] \le \Ex[e^{-r \tau \wedge \tau^x_{\xss (y)}}
  \phi(S^{x}_{\tau \wedge \tau^x_{\xss (y)} },y)]\le \underset{t \rightarrow +
    \infty}{\liminf} \Ex [e^{-r \tau \wedge \tau^x_{\xss (y)} \wedge t}
  \phi(S^{x}_{\tau \wedge \tau^x_{\xss (y)}\wedge t },y)]. \]
By Itô's formula,
  \begin{align*}
    de^{-rt}\phi(S^x_t,y)
    =&e^{-rt}\left[\frac{1}{2}\sigma^2(S^x_t)(S^x_t)^2\partial^2_x
    \phi(S^x_t,y) +(r-\delta)S^x_t\partial_x \phi(S^x_t,y) -r\phi(S^x_t,y) \right]dt \\
  &+e^{-rt}\partial_x \phi(S^x_t,y)\sigma(S^x_t)S^x_tdW_t
  \end{align*}

  On $\{t < \tau^x_{\xss (y)} \}$, we have $S^{x}_t<\xss(y)<X(y)$. In case $r\ge
  \delta$, it is obvious that the drift term is nonpositive since it is a sum of three
  nonpositive terms. This ensures that for $t\geq 0$, $\Ex [e^{-r \tau \wedge \tau^x_{\xss (y)} \wedge t}
  \phi(S^{x}_{\tau \wedge \tau^x_{\xss (y)}\wedge t },y)]\le \phi(x,y)$
  and $\tau_x^P=0$ is optimal.\\
In case $r \le \delta$, let us check that the
  drift term is still nonpositive by proving that the sum of the two last terms
  is nonpositive.\\The function
  is $x \mapsto
  (r-\delta)x\partial_x\phi(x,y)-r\phi(x,y)$ is nondecreasing on
  $(0,X(y))$ since its
  derivative is equal to $-\delta
  \partial_x\phi(x,y)+(r-\delta)x\partial^2_x\phi(x,y)$. Thus, using
  \eqref{EDf} for the second equality,  we have for $x \le x^*(y)$
  \begin{eqnarray*}
    (r-\delta)x\partial_x\phi(x,y)-r\phi(x,y) &\le
    &(r-\delta)\xss(y)\partial_x\phi(\xss(y),y)-r\phi(\xss(y),y) \\
    & =& \phi(\xss(y),y) \left[(r-\delta)\xss(y)\frac{f'(\xss(y))}{f(\xss(y))}-r\right]\\
    &=& -\frac{1}{2}\phi(\xss(y),y)\sigma^2(\xss(y))\xss(y)^2
    \frac{f''(\xss(y))}{f(\xss(y))} < 0.
  \end{eqnarray*}

\end{proof}

\begin{proposition}\label{opt_call}Let us fix a strike $x>0$ and assume $Y(x)>0$. If $Y(x)=+\infty$, then $\forall y>0,  c_\eta(y,x)=0$. Otherwise there is a unique $\yse(x) \in
  (Y(x),+\infty)$ such that $\tau^c_x=\inf \{t \ge 0, \bar{S}^y_t \ge \yse(x) \}$ (convention $\inf
 \emptyset=+\infty$) is an optimal stopping time for $c_\eta$ and:
 \begin{equation}\label{DEFy*}
  \forall y \ge \yse(x), c_\eta(y,x)=\phi(x,y), \ \forall y <
  \yse(x), c_\eta(y,x)= \frac{\phi(x,\yse(x)) }{
    g (\yse (x))}g (y) > \phi(x,y).  \hspace{0.1cm}
\end{equation}
In addition, we have $\frac{\phi (x,\yse(x))}{\partial_y \phi (x,\yse(x))}=
\frac{g(\yse(x))}{g'(\yse(x))}$ and the smooth-fit principle
holds.
Last, $x\in \{ z:0<Y(z)<\infty\} \mapsto \yse(x)$ is~$\mathcal{C}^1$ and
satisfies the following ODE:
\begin{eqnarray}\label{EDOy*}
  \yse(x)'&=& \left[\partial^2_{xy}\phi (x,\yse(x)) - \frac{\partial_x \phi
      (x,\yse(x))\partial_y \phi (x,\yse(x))}{\phi (x,\yse(x))}\right] \\
  && \times \frac{\yse(x)^2 \eta^2(\yse(x))}{2[\delta\phi (x,\yse(x)) +(r-\delta)\yse(x)
    \partial_y \phi (x,\yse(x)) ]-\yse(x)^2 \eta^2(\yse(x))\partial^2_y \phi
    (x,\yse(x))}. \nonumber \end{eqnarray}
It is strictly
increasing if one assume moreover $\phi \partial_x\partial_y \phi >\partial_x \phi\partial_y \phi $ on $\{\phi(x,y)>0\}$.
\end{proposition}
\begin{proof}
We introduce $h(y)=\frac{\phi(x,y)}{g(y)}$ which vanishes for
$y\le Y(x)$ and for $y=+\infty$. Indeed, the concavity ensures that $y \mapsto \phi(x,y)$ is
bounded from above by some linear function and we have already shown in \cite{AJ} that
$g(y)\underset{y \rightarrow + \infty}{\geq}cy^{1+a}$ for some $a,c>0$. We then obtain easily $\frac{\partial_y \phi (x,\yse(x))}{\phi (x,\yse(x))}=
\frac{g'(\yse(x))}{g(\yse(x))}$. The uniqueness of $\yse(x) \in
  (Y(x),+\infty)$ and the
  optimality of $\tau^c_x$ can be checked by arguments similar to the
  ones given in the proof of Proposition~\ref{opt_put}. To obtain the ODE~(\ref{EDOy*}) satisfied
  by~$\yse$, the calculations are the same as for~\eqref{EDOx*} in Proposition~\ref{opt_put}
exchanging $r \leftrightarrow \delta$, $\sigma \leftrightarrow \eta$ and $\partial_x
 \leftrightarrow \partial_y$.
\end{proof}
\begin{remark}\label{rem_1}
  We incidentally obtain in the proof of Proposition~\ref{opt_put} that \[\forall x
  \le \xss(y),  
  (r-\delta)x\partial_x\phi(x,y)-r\phi(x,y)<0 .\] Similarly,
  \[\forall y
  \ge \yse(x),  
  (\delta-r)y\partial_y\phi(x,y)-\delta\phi(x,y)<0. \]
  In particular, thanks to~\eqref{hypphi}, the
denominator of the second term in the r.h.s. of (\ref{EDOx*})
(resp. (\ref{EDOy*})) is positive.
\end{remark}

\subsection{Estimates on the exercise boundaries}

Now, we would like to get also estimations on the exercise boundaries. As in
\cite{AJ}, we use a comparison to the Black-Scholes model with constant
volatility for which estimations are
easier to get.

\begin{proposition}\label{compfront}
  Let us consider two volatility functions $\sigma_1$ and $\sigma_2$ (resp.$\eta_1$
  and $\eta_2$ ) satisfying $(\Hsw)$ such that   $\forall x>0, \sigma_1(x) \le  \sigma_2(x)$ (resp.
  $\forall x>0, \eta_1(x) \le  \eta_2(x)$). Then, we have:
\[  \forall x,y >0, P_{\sigma_1}(x,y) \le  P_{\sigma_2}(x,y) \text{ (resp. } \forall x,y >0,
c_{\eta_1}(y,x) \le  c_{\eta_2}(y,x)\ )\]
and we can compare the exercise boundaries:
\[ \forall y>0, x^*_{\sigma_2}(y) \le x^*_{\sigma_1}(y) \text{ (resp. } \forall x>0, y^*_{\eta_1}(x) \le  y^*_{\eta_2}(x) \ ). \]
\end{proposition}
\begin{proof}
Let us focus on the put case. If $P_{\sigma_1}(x,y)=\phi(x,y)$, we have
clearly $P_{\sigma_1}(x,y) \le  P_{\sigma_2}(x,y)$. Otherwise we have  $P_{\sigma_1}(x,y)=\phi(x^*_{\sigma_1}(y),y) \Ex[e^{-r
  \tau^{x,\sigma_1}_{x^*_{\sigma_1}(y)}}]$ where for $i\in\{1,2\}$,
$\tau^{x,\sigma_i}_{z}=\inf\{t\geq 0:S^{x,i}_t=z\}$ with $S^{x,i}_t$
solving \eqref{EDS} for the volatility function $\sigma_i$. Thanks to~(\ref{convex_fg}), we
  know that $f_{\sigma_1}$ is a convex function. According to the proof of Proposition~1.9~\cite{AJ}, $\Ex[e^{-r
  \tau^{x,\sigma_1}_{x^*_{\sigma_1}(y)}}] \le \Ex[e^{-r
    \tau^{x,\sigma_2}_{x^*_{\sigma_1}(y)}}]$. Therefore,
    \[ P_{\sigma_1}(x,y) \le \phi(x^*_{\sigma_1}(y),y) \Ex[e^{-r
    \tau^{x,\sigma_2}_{x^*_{\sigma_1}(y)}}] \le  P_{\sigma_2}(x,y). \]
\end{proof}

\begin{proposition}\label{bound_boundaries} Let $\overline{\sigma}$
  (resp. $\overline{\eta}$) denote an
   upper bound of the function $\sigma(.)$ (resp. $\eta(.)$). Then, \begin{align*}
   &\forall y>0\mbox{ s.t. }X(y)>0,\;\frac{a(\overline{\sigma})}{a(\overline{\sigma})-1} X(y) \le \xss(y)
   < X(y)\\
\bigg(\mbox{resp. }&\forall x>0\mbox{ s.t. }0<Y(x)<+\infty,\; 
   Y(x) < \yse(x) \le \frac{b(\overline{\eta})}{b(\overline{\eta})-1}Y(x)\bigg)
  \end{align*}
where  $a(\varsigma)=\frac{\delta-r+\varsigma^2/2 -
  \sqrt{(\delta-r+\varsigma^2/2)^2 + 2 r \varsigma^2}}{\varsigma^2}$
is an increasing function on $(0,+\infty)$ such that
$\lim_{\varsigma\rightarrow+\infty}a(\varsigma)=0$ and 
\begin{equation*}
   \lim_{\varsigma\rightarrow 0}a(\varsigma)=\begin{cases}
      -\frac{r}{\delta-r}\mbox{ if }\delta>r\\
-\infty\mbox{ otherwise}
   \end{cases}
\end{equation*}(resp. $b(\varsigma)=1-a(\varsigma)>1$).
\end{proposition}
\begin{proof}
  When $\delta=r$, the properties of $a(\varsigma)=\frac{1}{2} -
  \sqrt{\frac{1}{4}+
    \frac{2 r }{\varsigma^2}}$ are obvious. Otherwise,
  $a(\varsigma)=A(\frac{\delta-r}{\varsigma^2})$ with
  $A(x)=x+\frac{1}{2}-\sqrt{(x+\frac{1}{2})^2+\frac{2rx}{\delta-r}}$. Remarking that $\lim_{x\rightarrow -\infty}A(x)=-\infty$, $A(0)=0$ and
  $\lim_{x\rightarrow +\infty}A(x)=-\frac{r}{\delta-r}$, one
  easily deduces the limits of $a(\varsigma)$ as $\varsigma$ tends to
  $0$ or $+\infty$. Since
  $A'(x)=\frac{\sqrt{(x+\frac{1}{2})^2+\frac{2rx}{\delta-r}}-\left(x+\frac{1}{2}+\frac{r}{\delta-r}\right)}{\sqrt{(x+\frac{1}{2})^2+\frac{2rx}{\delta-r}}}$ and $(x+\frac{1}{2})^2+\frac{2rx}{\delta-r}-\left(x+\frac{1}{2}+\frac{r}{\delta-r}\right)^2=\frac{-r\delta}{(\delta-r)^2}\leq 0$, $A'(x)$ has the same sign as $-\left(x+\frac{\delta+r}{2(\delta-r)}\right)$. In particular $A'$ is negative on $(0,+\infty)$ when $\delta >r$ and positive on $(-\infty,0)$ when $\delta<r$. One easily deduces the monotonicity properties of $a$.

Let us deduce the estimation for the put case. Thanks to Proposition
\ref{compfront}, we have
  $\xss(y) \ge x^*_{\overline{\sigma}}(y)$. The solution of the
  EDO~(\ref{EDf}) with a volatility function constant equal to $\overline{\sigma}$ is
  $f_{\overline{\sigma}}(x)=x^{a(\overline{\sigma})}$. Let us consider the function
  $ x\in(0,X(y)) \mapsto \frac{\phi(x,y)}{\partial_x \phi(x,y)}- \frac{f_{\overline{\sigma}}(x)}{f'_{\overline{\sigma}}(x)}$. Its
  derivative $\frac{-\phi(x,y)\partial^2_x \phi(x,y)}{(\partial_x
    \phi(x,y))^2}+\frac{f_{\overline{\sigma}}(x)f''_{\overline{\sigma}}(x)}{(f'_{\overline{\sigma}}(x))^2}$ is greater than $\frac{a(\overline{\sigma})-1}{a(\overline{\sigma})}$ since $\partial_x^2\phi\leq 0$ and $\frac{f_{\overline{\sigma}}(x)f''_{\overline{\sigma}}(x)}{(f'_{\overline{\sigma}}(x))^2}=\frac{a(\overline{\sigma})-1}{a(\overline{\sigma})}$. Integrating this inequality between $x^*_{\overline{\sigma}}(y)$ and $X(y)$ then using \eqref{defx*} and remarking that by \eqref{hypphi}, $\partial_x\phi(X(y)^-,y)<0$ and
$\frac{\phi}{\partial_x\phi}(X(y)^-,y)=0$, we get
  $ - \frac{1}{a(\overline{\sigma})}X(y)  \ge
  \frac{a(\overline{\sigma})-1}{a(\overline{\sigma})}
  (X(y)-x^*_{\overline{\sigma}}(y))$ and thus: \[x^*_{\overline{\sigma}}(y)\ge \frac{a(\overline{\sigma})}{a(\overline{\sigma})-1} X(y) . \]
  The proof for $\yse$  works in the same way considering the function $y
  \mapsto  \frac{\phi(x,y)}{\partial_y\phi(x,y)}-
  \frac{g_{\overline{\eta}}(y)}{g'_{\overline{\eta}}(y)}$.\\

\end{proof}
\begin{remark}\label{surHvol} In the Call-Put case $\phi(x,y)=(y-x)^+$,
  since $\partial^2_x\phi(x,y)=0$ for $x<X(y)=y$, under $(\Hs)$ one
  obtains $\xss(y)\leq
  \frac{a(\underline{\sigma})}{a(\underline{\sigma})-1}y$ by an easy
  adaptation of the arguments given in the proof of
  Proposition~\ref{bound_boundaries}. In \cite{AJ}, this estimate
  combined with the ODE \eqref{EDOx*} derived below allowed us to characterise explicitly the set
of exercise boundaries $\xss$ and get a one-to-one correspondence
between the volatility
functions satisfying $(\Hs)$ and the exercise boundaries.\\
For general payoff functions $\phi$, because $\partial_x^2\phi$ does not
vanish, we were not able to get under $(\Hs)$ an upper-bound for $\xss$ better than
$\xss(y)<X(y)$ which already holds under $(\Hs')$. That is why we
work with hypothesis $(\Hs')$ in the present paper.\end{remark}

\section{Duality}\label{secduality}
Let us now investigate conditions ensuring
\begin{equation}
   \forall x,y>0,\;P_\sigma(x,y)=c_{\eta}(y,x)\label{dual}.
\end{equation} First, in order to use
the pricing formulas given in Propositions \ref{opt_put} and
\ref{opt_call}, we assume that for all $x>0$, $Y(x)>0$ condition which
implies $X(0^+)=0$.\\
Since $\Phi\neq \emptyset$, there exists $(x,y)\in\R_+^*\times \R_+^*$
such that $\phi(x,y)>0$. Then $X(y)>0$ and $Y(x)<+\infty$, and by
Propositions \ref{opt_put} and \ref{opt_call}, the functions $z\mapsto
P_\sigma(z,y)$ and $z\mapsto c_\eta(z,x)$ do not vanish on
$(0,+\infty)$. If for some $y'\in
(0,y)$, one had $X(y')=0$, then $\phi$
and therefore $P_\sigma$ would vanish on $(0,+\infty)\times (0,y']$. In
particular $P_\sigma$ would vanish on $\{x\}\times (0,y']$ preventing
\eqref{dual}. In the same way, if one had $X(+\infty)<+\infty$, then $c_\eta$ would vanish on
$(0,+\infty)\times[X(+\infty),+\infty)$ preventing \eqref{dual}.
That is why we make the following assumption on $X$:
\begin{equation}\label{Boundphi} \forall y>0, X(y)>0, \ X(0^+)=0 \text{ and } X(+\infty) =+\infty. \end{equation}
This assumption automatically ensures $Y(0^+)=0$, $0<Y(x)<+\infty$ for $x>0$ and
$Y(+\infty)=+\infty$.
We are now able give a necessary and sufficient condition for
\eqref{dual} to hold.
\begin{theorem}\label{thmdual}
  Assume that $\phi$ satisfy (\ref{Boundphi}) and that $\sigma$ and $\eta$ satisfy $(\Hs')$. Then,
  \eqref{dual} holds if and only if $\xss$ and $\yse$ are increasing
  reciprocal functions. 
\end{theorem}
\begin{proof} 
  This result can be
  checked by an immediate adaptation of the proof of Theorem~4.1~\cite{AJ}, except
  for the increasing property in the necessary condition that we explain here. The
  equality of the exercise regions writes
  \begin{equation}\label{egal_ex}
    \{(x,y)\in (\mathbb{R}^*_+)^2: x
    \le \xss(y)\}=\{(x,y)\in (\mathbb{R}^*_+)^2: \yse(x) \le y \},\end{equation}
  and thus $x\le \xss(
  \yse(x))$. Therefore $(x',\yse(x))$ belongs to the exercise region for $x' \le x$
  and we get $\yse(x')\le \yse(x)$. Similarly, $\xss$ is nondecreasing.
  Therefore, using Propositions \ref{opt_put},
  \ref{opt_call} and \ref{bound_boundaries} , we get that $\xss$ and $\yse$ are 
  continuous nondecreasing functions from $\mathbb{R}^*_+$ onto $\mathbb{R}^*_+$.
  From~\eqref{egal_ex}, they are reciprocal functions. Since they are both
  continuous, they are increasing.
\end{proof}

Let us recall here that  under the following assumption on  $\phi$ 
\begin{equation}\label{bound_incr} \phi \partial^2_{xy} \phi >\partial_x \phi\partial_y \phi \text{ on } \Phi,  \end{equation}
Propositions~\ref{opt_put} and~\ref{opt_call} ensure that the exercise boundaries are
automatically increasing. We give a general class of functions $\phi$ that satisfy all the
required assumptions.
\begin{example}
   Let $\psi:\mathbb{R}_+ \rightarrow \mathbb{R}_+$ be
  an increasing concave function $\mathcal{C}^2$ on $(0,+\infty)$ and
  such that $\psi(0)=\psi(0^+)=0$, $\psi_x:\R_+^*\rightarrow \R_+^*$
  (resp. $\psi_y:\R_+^*\rightarrow \R_+^*$) be a $C^2$ increasing convex
  (resp. concave) function such that $\psi_x(0^+)=0$
  (resp. $\psi_y(0^+)=0$ and $\psi_y(+\infty)=+\infty$). Then the function
  $\phi(x,y)=\psi((\psi_y(y)-\psi_x(x))^+)$ satisfies~\eqref{hypphi}. It is such that
  $X(y)=\psi_x^{-1}(\psi_y(y))$, $Y(x)=\psi_y^{-1}(\psi_x(x))$ and (\ref{Boundphi}) and (\ref{bound_incr}) hold.
\end{example}

For some specific payoff functions of this family, we are now going to state conditions
on $\sigma$ and $\eta$ such that $\xss$ and $\yse$ are reciprocal functions. We
first recall results obtained in \cite{AJ} in the call-put case
$\phi(x,y)=(y-x)^+$. Then, we address two generalizations :
$\phi(x,y)=(\psi_y(y)-\psi_x(x))^+$ and $\phi(x,y)={(y-x)^+}^\gamma$.

\subsection{The call-put case $\phi(x,y)=(y-x)^+$}\label{callput}

Let us recall here the main result obtained in~\cite{AJ}:
\begin{theorem}\label{DualProp_cp}
 Let us consider two volatility functions satisfying $(\Hs)$. The following conditions are equivalent:\\
 (1)  $   \forall x,y>0,\  P_\sigma(x,y)= c_\eta(y,x).$\\
(2) $\eta\equiv\tsigma$ where $\tsigma(y)= 2(y-\xss(y))(ry-\delta \xss(y))/[y \xss(y)
    \sigma (\xss(y))].$\\
(3) $\sigma\equiv\uteta$ where $ \uteta(x)=
    2(\yse(x)-x)(r \yse(x)-\delta x) / [ \yse(x) x
    \eta (\yse(x))].$
\end{theorem}
As proved in \cite{AJ}, if $\sigma$ (resp. $\eta$)
satisfies  $(\Hs)$, then $\tsigma$ (resp. $\uteta$) also satisfies
$(\Hs)$. This very convenient property ensures that for a given volatility
function $\sigma$ satisfying $(\Hs)$, there always exists a dual
volatility function $\eta$ also satisfying $(\Hs)$ such that
condition~(1) above holds. Unfortunately, for the more general payoff
functions that we consider in the sequel, stability of the
Hypotheses $(\Hs)$ or $(\Hsw)$ is no longer straightforward. And it may
happen that no dual volatility function $\eta$ can be associated with $\sigma$.

\subsection{The case $\phi(x,y)=(\psi_y(y)-\psi_x(x))^+$}

In this section, we will focus on the case $\phi(x,y)=(\psi_y(y)-\psi_x(x))^+$ where $
\psi_y:\mathbb{R}^*_+ \rightarrow \mathbb{R}^*_+$ (resp. $\psi_x:\mathbb{R}^*_+ \rightarrow
\mathbb{R}^*_+$) is a $C^2$ increasing concave (resp. convex)
function such that $\psi_y(0^+)=0$ and $\psi_y(+\infty)=+\infty$ (resp. $\psi_x(0^+)=0$). Then one has $X(y)=\psi_x^{-1}(\psi_y(y))$ and $Y(x)=\psi_y^{-1}(\psi_x(x))$.
 
Let us first give an example of application of Theorem \ref{thmdual}
when $\psi_x$ and $\psi_y$ are power functions and the local volatility
functions $\sigma$ and $\eta$ are constant.
\begin{example}\label{Ex_BS}
   Let us suppose that $\phi(x,y)=(y^{\gamma'}-x^{\gamma})^+$ where
   $\gamma' \in (0,1]$ and $\gamma \geq 1$. When the local volatility
   function $\sigma$ is a constant and equal to $\varsigma$,
   $f(x)=x^{a(\varsigma)}$ with $a(\varsigma)$ given in Proposition
   \ref{bound_boundaries}. The equality $\frac{\partial_x \phi (\xss(y),y)}{\phi (\xss(y),y)}=
\frac{f'(\xss(y))}{f(\xss(y))}$ then yields
$\xss(y)=\left(\frac{a(\varsigma)}{a(\varsigma)-\gamma}\right)^{1/\gamma}y^{\gamma'/\gamma}$.
In the same way, for $\eta$ constant equal to $\nu$, as
$g(x)=x^{b(\nu)}$ with $b(\nu)=1-a(\nu)$,
$\yse(x)=\left(\frac{b(\nu)}{b(\nu)-\gamma'}\right)^{1/\gamma'}x^{\gamma/\gamma'}$.
These boundaries are reciprocal functions as soon as 
\begin{equation*}
   \gamma' a(\varsigma)+\gamma b(\nu)=\gamma\gamma'.
\end{equation*}
According to Proposition \ref{bound_boundaries}, when
$r\geq \delta$, for fixed $\varsigma\in (0,+\infty)$ this equation
admits a solution $\nu\in (0,+\infty)$ iff
$\varsigma<a^{-1}(\gamma(1-\frac{1}{\gamma'}))$ and it admits a
solution $\varsigma\in (0,+\infty)$ for any fixed $\nu\in
(0,+\infty)$. When $\delta>r$, there is no solution if
$\gamma(1-\frac{1}{\gamma'})\leq -\frac{r}{\delta -r}$ and otherwise it
admits a solution $\nu$ for fixed
$\varsigma<a^{-1}(\gamma(1-\frac{1}{\gamma'}))$ and a solution
$\varsigma$ for fixed $\nu>b^{-1}(\gamma'(1+\frac{r}{\gamma(\delta-r)}))$.
\end{example}

For general functions $\psi_x$ and $\psi_y$, we are able to investigate
uniqueness for the ODEs ~(\ref{EDOx*}) and~(\ref{EDOy*}) which respectively write :
\begin{eqnarray}\label{EDOx*_2}
  \xss(y)'&=& \frac{\psi_x'(\xss(y)) \psi_y'(y)}{ \psi_y(y)-\psi_x(\xss(y))} \\&& \times \frac{\xss(y)^2 \sigma^2(\xss(y))}{2[r (\psi_y(y)-\psi_x(\xss(y))) +(r-\delta)\xss(y)
   \psi_x'(\xss(y)) ]+\xss(y)^2 \sigma^2(\xss(y)) \psi_x''
   (\xss(y))} \nonumber\\
\label{EDOy*_2}
  \yse(x)'&=& \frac{\psi_x'(x) \psi_y'(\yse(x))}{\psi_y(\yse(x))-\psi_x(x)} \\ && \times \frac{\yse(x)^2 \eta^2(\yse(x))}{2[\delta (\psi_y(\yse(x))-\psi_x(x)) +(r-\delta)\yse(x)
     \psi_y' (\yse(x)) ]-\yse(x)^2 \eta^2(\yse(x)) \psi_y''
    (\yse(x))} \nonumber
.\end{eqnarray}
\begin{proposition}\label{uniy*}When $\eta$ satisfies $(\Hsw)$, the boundary $\yse(x)$ is the unique solution $y(x)$
  of~(\ref{EDOy*_2}) on $\mathbb{R}^*_+$ that is
  increasing, such that $Y(x)<y(x)$ and $y(0^+)=0$.
\end{proposition}
\begin{proof}
  Let us consider $y_1(x)$ and $y_2(x)$, two solutions of \eqref{EDOy*_2} on $\mathbb{R}^*_+$ that are
  increasing and such that $y_i(x)>Y(x)$ and $y_i(0^+)=0$ for $i\in \{1,2\}$. In
  particular, $y_1$ and $y_2$ are bijections on $\mathbb{R}^*_+$ and we
  may define $\tilde{I}(x)=y_1^{-1}(y_2(x))/x$. By an easy computation,
  one checks
 \begin{align*}
  &\tilde{I}'(x)= \frac{1}{x} \Bigg(
     \frac{\psi_x'(x)}{\psi_x'(x\tilde{I}(x))}\times\left[1+\frac{\psi_x(x)-\psi_x(x\tilde{I}(x))}{\psi_y(y_2(x))-\psi_x(x)}\right] \\& \times \left[1+\frac{2\delta [\psi_x(x)-\psi_x(x\tilde{I}(x))]}{2[\delta (\psi_y(y_2(x))-\psi_x(x)) +(r-\delta)y_2(x)
     \psi_y' (y_2(x)) ]-y_2(x)^2 \eta^2(y_2(x)) \psi_y''
     (y_2(x))}\right]-\tilde{I}(x)\Bigg).
\end{align*}
The constant $1$ is clearly solution to this equation and we want
to check that $\tilde{I} \equiv 1$. Let us suppose that $\tilde{I} \not \equiv 1$. Thanks to the Cauchy-Lipschitz
theorem, it induces that either $\forall x>0, \tilde{I}(x)>1$ or $\forall x>0, \tilde{I}(x)<1$. Let us
suppose $\forall x>0, \tilde{I}(x)<1$. Then, it is easy to see from the last expression
that
\[ \forall x>0, \ \tilde{I}'(x) \ge \frac{1}{x}(1 - \tilde{I}(x) ). \]
Indeed, since $\psi_x'$ is non decreasing, we have
$\psi_x'(x) \ge \psi_x'(x\tilde{I}(x))$ and the terms into brackets are also greater than
$1$ because $\psi_x$ is increasing and both denominators are nonnegative
as  $\psi_y(y_2(x))>\psi_x(x)$, $y_2'(x)\ge 0$ and $y_2$ solves \eqref{EDOy*_2}.
In particular we have shown that $\tilde{I}'(x)  >0$ and therefore,
\[ \forall x\in (0,1), \ \tilde{I}'(x) \ge \frac{1}{x}(1 - \tilde{I}(1) ). \]
Thus, we get $\tilde{I}(1)-\tilde{I}(x) \ge ( \tilde{I}(1) -1 )\ln(x)
\underset{x\rightarrow 0^+}{\rightarrow}+\infty$ and so $\tilde{I}(x)\underset{x\rightarrow
  0^+}{\rightarrow}-\infty$ which is contradictory since $\forall x>0,
\tilde{I}(x)>0$.

When $\forall x>0, \tilde{I}(x)>1$, considering
$y_2^{-1}(y_1(x))/x$ instead of $\tilde{I}(x)$, we get the same contradiction as previously. \end{proof}

Now let us turn to the uniqueness result on the boundary $\xss(y)$.
\begin{proposition}\label{unix*}
  Let $\sigma$ satisfy $(\Hsw)$ and $\psi_x$ be such that :
  \begin{equation}\label{Hyp_psix}
    \forall \alpha \in (0,1), \exists C_\alpha>0, \forall x>0, \  \psi_x(\alpha x) \ge
    C_\alpha \psi_x(x).
  \end{equation}
  
  The boundary $\xss(y)$ is the unique solution $x(y)$
  of~(\ref{EDOx*_2}) on $\R_+^*$ that is
  increasing and such that $ \exists \alpha \in(0,1)$,
 \begin{equation}
   \forall y>0, \  \alpha X(y) \le x(y) < X(y).\label{bornfront}
 \end{equation}
\end{proposition}
Hypothesis \eqref{Hyp_psix} is satisfied by the function $x^a$ with $a
\ge 1$ but not by the function $\exp(bx)-1$ with $b>0$.
\begin{proof}The boundary $\xss(y)$ satisfies \eqref{bornfront} with
  $\alpha=\frac{a(\overline{\sigma})}{a(\overline{\sigma})-1}$
according to Proposition \ref{bound_boundaries}.\\
  Let $x_1$ and $x_2$ denote two solutions of~(\ref{EDOx*_2}) satisfying
  \eqref{bornfront} with respective constants $\alpha_1,\alpha_2\in
  (0,1)$ and
  $\hat{I}(y)=\psi_y(x_1^{-1}(x_2(y)))/\psi_y(y)$. One has
\begin{align*}
 &\hat{I}'(y)=\frac{\psi_y'(y)}{\psi_y(y)} \times \Bigg(\left[ 
    \frac{\hat{I}(y)-\psi_x(x_2(y))/\psi_y(y)}{1-\psi_x(x_2(y))/\psi_y(y)} \right] \\
  &\times\left[1+ \frac{2r \psi_y(y)(\hat{I}(y)-1) }{2[r (\psi_y(y)-\psi_x(x_2(y))) +(r-\delta)x_2(y)
   \psi_x'(x_2(y)) ]+x_2(y)^2 \sigma^2(x_2(y)) \psi_x''
   (x_2(y))}\right]- \hat{I}(y) \Bigg).
\end{align*}
Let us suppose that $\hat{I}(y) \not \equiv 1$. Thanks to the Cauchy-Lipschitz
theorem, we have either $\forall y>0,  \hat{I}(y)>1$ or $\forall y>0,  \hat{I}(y)<1$.
Let us suppose that $\forall y>0,  \hat{I}(y)>1$. As in the last proof, the second bracket is
greater than $1$, and we get
\[\forall y>0, \ \hat{I}'(y) \ge \frac{\psi_y'(y)}{\psi_y(y)}(\hat{I}(y)-1)
\frac{\psi_x(x_2(y))/\psi_y(y)}{1-\psi_x(x_2(y))/\psi_y(y)}. \]
Since $x_2(y)<X(y)=(\psi_x)^{-1}(\psi_y(y))$, we have $0<\psi_x(x_2(y))/\psi_y(y)<1$
and therefore $\hat{I}'(y)>0$.
Since $x_2$ satisfies \eqref{bornfront} with constant $\alpha_2$, we get by \eqref{Hyp_psix}
\[\forall y>0,\;\psi_x(x_2(y)) \ge \psi_x(\alpha_2 X(y)) \ge C_{\alpha_2} \psi_y(y) .\]
Since $z\mapsto \frac{z}{1-z}=-1+\frac{1}{1-z}$ is increasing on $(0,1)$ and
$\hat{I}$ is increasing, we deduce that
\[\forall y \ge 1, \ \hat{I}'(y) \ge \frac{\psi_y'(y)}{\psi_y(y)}(\hat{I}(1)-1)
\frac{C_{\alpha_2}}{1-C_{\alpha_2}}. \]
As a consequence,
\[\hat{I}(y)-\hat{I}(1)\ge (\hat{I}(1)-1)
\frac{C_{\alpha_2}}{1-C_{\alpha_2}} \ln\left(\frac{\psi_y(y)}{\psi_y(1)}\right) \underset{y \rightarrow
  + \infty}{\rightarrow} + \infty.\]
In the same time, since $x_1$ satisfies \eqref{bornfront} with constant
$\alpha_1$, we have $ X(y)=(\psi_x)^{-1}(\psi_y(y))  \le
\frac{1}{\alpha_1}x_1(y)$ and therefore $x_1^{-1}(x) \le (\psi_y)^{-1}(\psi_x(\frac{x}{\alpha_1}))$. We get
$\psi_y(x_1^{-1}(x_2(y)))\le \psi_x(\frac{x_2(y)}{\alpha_1} )\le \frac{\psi_x( x_2(y))}{C_{\alpha_1}}<\frac{\psi_y(y)}{C_{\alpha_1}} $ and
thus $\hat{I}(y) \le \frac{1}{C_{\alpha_1}}$ which is contradictory with $\hat{I}(+ \infty)=+ \infty.$

When $\forall y>0,  \hat{I}(y)<1$, considering
$\psi_y(x_2^{-1}(x_1(y)))/\psi_y(y)$ instead of $\hat{I}(y)$, we
reach the same contradiction as previously.
\end{proof}

Like in the call-put case, we are now able to state a more precise duality result.

\begin{theorem}\label{dualpart}
Let us assume  that
$\sigma$ and $\eta$ satisfy $(\Hsw)$ and set\\
$A(y)= \left[ \frac{ \psi_y(y)-\psi_x(\xss(y))} {\psi_x'(\xss(y)) \psi_y'(y)}\right]^2 \times \frac{2[r (\psi_y(y)-\psi_x(\xss(y))) +(r-\delta)\xss(y)
   \psi_x'(\xss(y)) ]+\xss(y)^2 \sigma^2(\xss(y)) \psi_x''
   (\xss(y))}{\xss(y)^2 \sigma^2(\xss(y))}$,\\
$B(x)= \left[\frac{\psi_y(\yse(x))-\psi_x(x)}{\psi_x'(x) \psi_y'(\yse(x))}\right]^2 \times \frac{2[\delta (\psi_y(\yse(x))-\psi_x(x)) +(r-\delta)\yse(x)
     \psi_y' (\yse(x)) ]-\yse(x)^2 \eta^2(\yse(x)) \psi_y''
    (\yse(x))}{\yse(x)^2 \eta^2(\yse(x))}$ which are positive functions
  according to Remark~\ref{rem_1}.
Then, the following assertions are equivalent~:\\
(1) $\forall
  x,y>0, \ P_{\sigma}(x,y)=c_\eta(y,x)$.\\
(2) $\forall y>0$, $\min[1+\psi''_y(y)A(y),\delta (\psi_y(y)-\psi_x(\xss(y))) +(r-\delta)y
     \psi_y' (y)]>0$ and  $\eta \equiv \tsigma$ where
  \begin{equation}\label{sigma_tilde}
    \tsigma(y)=\frac{1}{y}\sqrt{2[\delta (\psi_y(y)-\psi_x(\xss(y))) +(r-\delta)y
     \psi_y' (y) ]\frac{A(y)}{1+\psi_y''(y)A(y)}}.
  \end{equation}
If one assumes moreover that $\psi_x$ satisfies~(\ref{Hyp_psix}), they are
also equivalent to\\
(3) $\forall x>0$,
  $\min[1-\psi_x''(x)B(x),r (\psi_y(\yse(x))-\psi_x(x)) +(r-\delta)x
   \psi_x'(x) ]>0$ and $\sigma \equiv \uteta$ where
  \begin{equation}\label{eta_utilde}
    \uteta(x) = \frac{1}{x} \sqrt{2[r (\psi_y(\yse(x))-\psi_x(x)) +(r-\delta)x
   \psi_x'(x) ] \frac{B(x)}{1-\psi_x''(x)B(x)} }.
  \end{equation}
\end{theorem}
Notice that one easily recovers the call-put formulas given in Theorem \ref{DualProp_cp} if one takes
$\psi_x(x)=x$ and $\psi_y(y)=y$. 
\begin{proof}
  Since the payoff function satisfies~\eqref{bound_incr}, by Theorem~\ref{thmdual}
  the assertion (1) is equivalent to the 
  reciprocity of the functions $\xss$ and $\yse$. Therefore the
  implications $(1)\Rightarrow(2)$ and $(1) \Rightarrow (3)$ are
  obtained by combining respectively 
  $(\yse)'(\xss(y))\xss(y)'=1$ and $(\xss)'(\yse(x))\yse(x)'=1$ with  (\ref{EDOx*_2})
  and (\ref{EDOy*_2}), the positivity of the terms between brackets in
  \eqref{sigma_tilde} and \eqref{eta_utilde} coming from Remark
  \ref{rem_1}.\\
Let us prove $(2)\Rightarrow(1)$. Computing $({\xss}^{-1})'(x)$
thanks to (\ref{EDOx*_2}), then using \eqref{sigma_tilde} written at the
point $y= {\xss}^{-1}(x)$, we check that ${\xss}^{-1}$
solves the same ODE as $\yse$. Since ${\xss}^{-1}$ is increasing and
${\xss}^{-1}(0^+)=0$, we conclude by
Proposition~\ref{uniy*} that ${\xss}^{-1}\equiv\yse$.

To prove $(3)\Rightarrow(1)$, we check in the same manner that
${\yse}^{-1}(y)$ solves the same ODE as $\xss(y)$. The function
$\yse(x)$ is increasing and according
to Proposition~\ref{bound_boundaries}, $\yse(x)\le \beta Y(x)$ for
$\beta=\frac{b(\overline{\eta})}{b(\overline{\eta})-1}>1$. With the
concavity of $\psi_y$ and $\psi_x^{-1}$ combined with
$\psi_y(0^+)=\psi_x^{-1}(0^+)=0$, this ensures\[\forall y>0,\;(\yse)^{-1}(y)\ge X(y/\beta)=(\psi_x)^{-1}(\psi_y(y/\beta))\ge
(\psi_x)^{-1}(\psi_y(y)/\beta)\ge X(y)/\beta.\] By
 Proposition~\ref{unix*}, we conclude that $(\yse)^{-1}\equiv\xss$. 
\end{proof}

To give an analytical example of non constant dual volatility functions,
we now assume that 
$\boxed{\phi(x,y)=(\alpha y- x^\gamma)^+}$ with $\alpha>0$ and $\gamma
  \ge 1$. For $a,b,c>0$, we introduce the reciprocal functions 
\[y^*(x)=\frac{1}{\alpha}x^\gamma \frac{x^\gamma+a}{bx^\gamma+c}\;\mbox{
  and }\; x^*(y)=\left[ \frac{1}{2} \left( b\alpha y - a +\sqrt{(b\alpha y - a)^2+4c \alpha y} \right) \right]^{1/\gamma}. \]
Under some assumptions on the coefficients $a,\ b$ and $c$, these
functions are the exercise boundaries associated with explicit dual
volatility functions.
\begin{proposition}
Let us assume that either $r\geq \delta$ and $\max(c/a,b) \le 1$ with
$\min(c/a,b) < 1 $ or $r < \delta$ and  $\max(c/a,b) \le \frac{1}{1+(\delta/r-1) \gamma}$ with
$\min(c/a,b) < \frac{1}{1+(\delta/r-1) \gamma}$. Let us also assume  $(\gamma -1) b(2c-a)+c(\gamma+1) \ge
0$. Then, the volatility functions
\begin{eqnarray*}
    \sigma(x) &=& \frac{1}{x} \sqrt{2[r (\alpha y^*(x)-x^\gamma) +(r-\delta)\gamma x^\gamma 
    ] \frac{B(x)}{1-\gamma(\gamma-1)x^{\gamma-2}B(x)} } \\  
&& \text{ with } B(x)=
\frac{1}{y^*(x)'}\left[\frac{\alpha y^*(x)-x^\gamma}{\alpha \gamma x^{\gamma-1}}\right],  \text{ and }\\
\eta(y)&=&\frac{1}{y}  \sqrt{2[ r \alpha y -\delta  x^*(y)^\gamma ] A(y)} \text{
  with } A(y)= \frac{1}{x^*(y)'}\left[ \frac{ \alpha y- x^*(y)^\gamma} { \alpha \gamma x^*(y)^{\gamma-1} }
\right]
\end{eqnarray*}
are well defined and satisfy $(\Hsw)$. Moreover, we have $\yse\equiv y^*$ and
$\xss\equiv x^*$ and thus the duality holds: $\forall x,y>0, \
P_\sigma(x,y)=c_\eta(y,x)$. 
\end{proposition}
When $r \ge \delta$, it is easy
to fulfill the required assumptions by taking for example $a$, $b$ and $c$ such that
$b<1$ and $1/2\le c/a<1$. When $r<\delta$, the first condition is
satisfied if $\max(c/a,b) < \frac{1}{1+(\delta/r-1) \gamma}$ and the second condition
can be rewritten $2\ge \frac{a}{c}-\frac{1}{b}\frac{\gamma+1}{\gamma-1}$. Thus taking
for example 
$b< \frac{1}{1+(\delta/r-1) \gamma}$  and then $\frac{c}{a}=b\frac{\gamma-1}{\gamma+1}$, one can get dual
volatility functions.

\begin{proof}
{\it First step: let us check that the functions $\sigma$ and $\eta$ are well defined and satisfy $(\Hsw)$.} Since
we have $\max(c/a,b) \le 1$ and $\min(c/a,b) < 1$, we get $y^*(x)
>\frac{1}{\alpha}x^\gamma=Y(x)$ (and thus $x^*(y) < (\alpha y) ^{1/\gamma} = X(y)$).
Since $y^*(x)'>0$ (and thus $x^*(y)'=1/{y^*}'(x^*(y))>0$), this ensures $B(x)>0$ and
$A(y)>0$. For $r \ge \delta$, it is then clear that  $r (\alpha
y^*(x)-x^\gamma) +(r-\delta)\gamma x^\gamma>0$ and $ r \alpha y -\delta
x^*(y)^\gamma >0$. For $\delta > r$, the condition $\max(c/a,b) \le \frac{1}{1+(\delta/r-1) \gamma}$
and $\min(c/a,b) < \frac{1}{1+(\delta/r-1) \gamma}$ ensures that $r (\alpha
y^*(x)-x^\gamma) +(r-\delta)\gamma x^\gamma>0$, but also $ r \alpha y -\delta
x^*(y)^\gamma >0$ (or equivalently $ r \alpha y^*(x) -\delta
x^\gamma >0$) since $ \frac{1}{1+(\delta/r-1) \gamma} \le r/\delta$ for $\gamma
\ge 1$. Thus, $\eta$ is well defined and positive. 
Since \[y^*(x)' =\frac{1}{\alpha} \gamma x^{\gamma-1} \frac{b x^{2\gamma} +2cx^{\gamma} +ac } {(bx^{\gamma}+c)^2}, \]
we get after some calculations that $\frac{B(x)}{1-\gamma(\gamma-1)x^{\gamma-2}B(x)}$
is equal to
\[  \frac{1}{\gamma} x^{2-\gamma}
\frac{(bx^\gamma+c)((1-b)x^\gamma+a-c)}{b(1+(\gamma-1)b)x^{2 \gamma}+((\gamma -1)
  b(2c-a)+c(\gamma+1) )x^\gamma+c(\gamma c +a-c) }\]
and is positive because we have assumed $(\gamma -1) b(2c-a)+c(\gamma+1) \ge
0$ (all other terms are positive). Thus $\sigma$ is well defined and we have


\begin{equation}
   \sigma(x)=\sqrt{\frac{2}{\gamma}   \frac{\left[r (\frac{x^\gamma+a}{bx^\gamma+c}-1) +(r-\delta)\gamma 
    \right](bx^\gamma+c)((1-b)x^\gamma+a-c)}{b(1+(\gamma-1)b)x^{2 \gamma}+((\gamma -1)
  b(2c-a)+c(\gamma+1) )x^\gamma+c(\gamma c +a-c) } }.\label{siggma}
\end{equation}
that is clearly bounded from above. To see that $\eta$ is also bounded from above we
calculate

\[ \eta(y)= \sqrt{2 \frac{r\alpha y - \delta x^*(y)^\gamma}{y} \times
  \frac{\alpha y - x^*(y)^\gamma}{\alpha^2 y}\times \frac{b x^*(y)^{2\gamma} +2cx^*(y)^{\gamma} +ac}{(bx^*(y)^\gamma+c)^2}}\]
using that $1/{x^*}'(y)={y^*}'(x^*(y)).$

{\it Second step:} We have $\frac{1}{\alpha} x^\gamma <y^*(x)\le \max(\frac{1}{b},\frac{a}{c})\frac{1}{\alpha} x^\gamma$
and thus $\max(\frac{1}{b},\frac{a}{c})^{-1/\gamma } (\alpha y) ^{1/\gamma} \le
x^*(y)< (\alpha y) ^{1/\gamma}$. From the definition of $\sigma$, we get $B(x)=\frac{x^2\sigma^2(x)}{2[r (\alpha y^*(x)-x^\gamma) +(r-\delta)\gamma x^\gamma 
  ]+ \gamma(\gamma-1)x^{\gamma}\sigma^2(x) }$. Combining this equality
for $x=x^*(y)$ with the definition of $B$, we deduce that $x^*$ solves
the ODE (\ref{EDOx*_2}). In the same manner, we show that $y^*$ solves the ODE
(\ref{EDOy*_2}). Thanks to Propositions~\ref{uniy*} and~\ref{unix*}, we conclude that $y^*\equiv\yse$ and  $x^*\equiv\xss$.
\end{proof}
\begin{remark}
  For $b=1$, we get cases where $\sigma$ and $\eta$ satisfy $(\Hsw)$ but not
  $(\Hs)$ since
$ \eta(y^*(x))=\sqrt{\frac{2}{\alpha} \frac{ (r-\delta b)x^\gamma + ra - \delta c
  }{x^\gamma+a} \times \frac{(b x^{2\gamma} +2cx^{\gamma}
    +ac)((1-b)x^\gamma+a-c)}{(x^\gamma+a)(bx^\gamma+c)^2}}$.
For $\delta>r$, $b=\frac{1}{1+(\delta/r-1) \gamma}$ and $\gamma>1$ we get cases where
$\eta$ satisfies $(\Hs)$ and $\sigma$ satisfies $(\Hsw)$ but not
  $(\Hs)$ (see \eqref{siggma}). If we have $\max(c/a,b)
< 1$ when $r\ge  \delta$ or $\max(c/a,b)<\frac{1}{1+(\delta/r-1) \gamma}$ when $r<
\delta$, one can check that $\sigma$ and $\eta$ satisfy $(\Hs)$.
\end{remark}

\begin{figure}[h]
  \psfrag{cas1}{\small \hspace{-2.5cm} \boxed{\phi(x,y)=((y-x)^+)^{\frac{3}{4}}}}
  \psfrag{cas2}{\small \hspace{-2.5cm} \boxed{\phi(x,y)=(0.97y-x^4)^+} }

  \centerline{\psfig{file=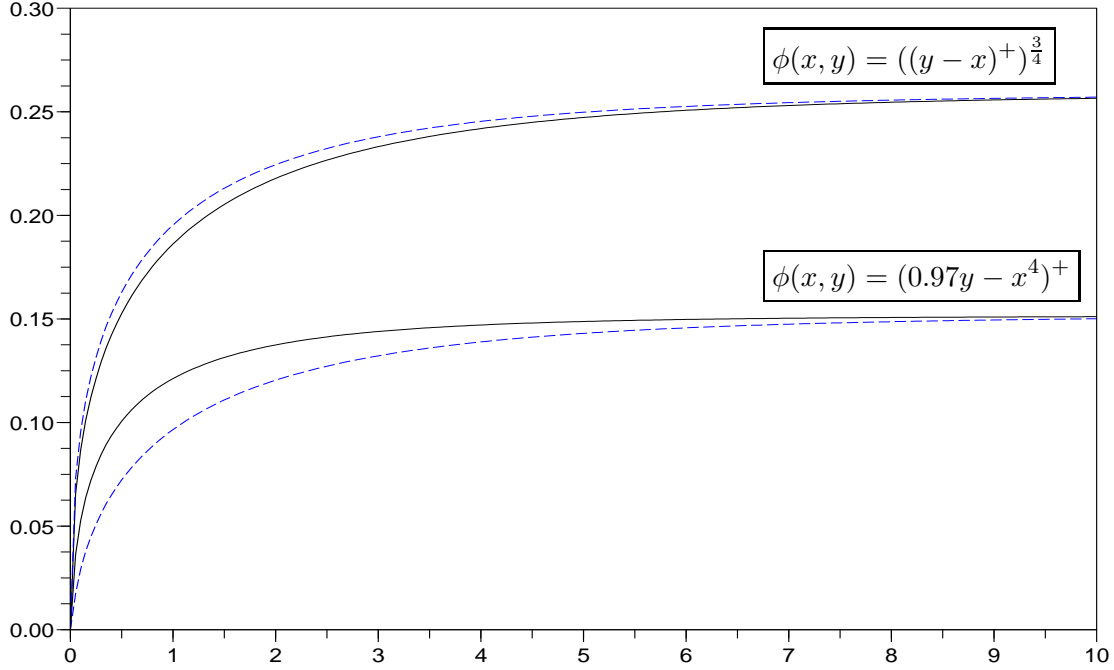,width=18cm,height=11cm} }
\caption{{\small $P_\sigma(T,x,y)$ (solid line) and $c_{\tsigma}(T,y,x)$ (dashed line) in function of the time
    $T$ for $a=1.5$, $b=5/9$, $c=1$, $x=1$, $y=0.99$, $r=0.2$ and $\delta=0.1$.}} \label{Dualgenfig}
\end{figure}

We have plotted in Figure~\ref{Dualgenfig} an example that illustrates
the duality. We have computed prices of American options with finite maturity $T$,
precisely $\underset{\tau \in
  \mathcal{T}_{0,T}} {\sup} \Ex \left[e^{-r\tau} \phi(S^{x}_{\tau},y)^+
\right]$ and $c_{\tsigma} (T,y,x)= \underset{\tau \in
   \mathcal{T}_{0,T}} {\sup} \Ex \left[e^{- \delta \tau}\phi(x,\overline{S}^{y}_{\tau})
 \right]$ where the supremum is taken over $\mathcal{T}_{0,T}$, the set of stopping
 times almost surely smaller than $T$. We
see that both converge to the same limit when $T$ is large.

\subsection{The case $\phi(x,y)={(y-x)^+}^\gamma$, $\gamma\in(0,1]$}
Let us first give an example of application of Theorem \ref{thmdual}
for this payoff when the local volatility
functions $\sigma$ and $\eta$ are constant.
\begin{example}\label{volconst}
   When the local volatility
   function $\sigma$ is a constant and equal to $\varsigma$,
   $f(x)=x^{a(\varsigma)}$ with $a(\varsigma)$ given in Proposition
   \ref{bound_boundaries}. The equality $\frac{\partial_x \phi (\xss(y),y)}{\phi (\xss(y),y)}=
\frac{f'(\xss(y))}{f(\xss(y))}$ then yields
$\xss(y)=\frac{a(\varsigma)}{a(\varsigma)-\gamma}y$.
In the same way, for $\eta$ constant equal to $\nu$, as
$g(x)=x^{b(\nu)}$ with $b(\nu)=1-a(\nu)$,
$\yse(x)=\frac{b(\nu)}{b(\nu)-\gamma}x$.
These boundaries are reciprocal functions as soon as 
\begin{equation*}
   a(\varsigma)+b(\nu)=\gamma.
\end{equation*}
According to Proposition \ref{bound_boundaries}, when
$r\geq \delta$, for fixed $\varsigma\in (0,+\infty)$ this equation
admits a solution $\nu\in (0,+\infty)$ iff
$\varsigma<a^{-1}(\gamma-1)$ and it admits a
solution $\varsigma\in (0,+\infty)$ for any fixed $\nu\in
(0,+\infty)$. When $\delta>r$, there is no solution if
$\gamma-1\leq -\frac{r}{\delta -r}$ and otherwise it
admits a solution $\nu$ for fixed
$\varsigma<a^{-1}(\gamma-1)$ and a solution
$\varsigma$ for fixed $\nu>b^{-1}(\gamma+\frac{r}{\delta-r})$.
\end{example}
For the particular choice $\phi(x,y)={(y-x)^+}^\gamma$ the ODEs \eqref{EDOx*}
and \eqref{EDOy*} write
\begin{align}
   \xss(y)'&=\frac{\gamma\xss(y)^2 \sigma^2(\xss(y))}{2[r(y-\xss(y))^2 +\gamma(r-\delta)\xss(y)
     (y-\xss(y)) ]+\gamma(1-\gamma)\xss(y)^2 \sigma^2(\xss(y))} \label{EDOx*_3}\\
\yse(x)'&=\frac{\gamma\yse(x)^2 \eta^2(\yse(x))}{2[\delta (\yse(x)-x)^2 +\gamma(r-\delta)\yse(x)
     (\yse(x)-x) ]+\gamma(1-\gamma)\yse(x)^2 \eta^2(\yse(x))}. \label{EDOy*_3}
\end{align}
Since $Y(x)=x$, $\yse(x)>x$ and by Remark~\ref{rem_1},
  $\delta(\yse(x)-x)+\gamma(r-\delta)\yse(x)>0$. It turns out that
uniqueness holds for the ODE~\eqref{EDOy*_3} under these conditions.
\begin{proposition}\label{uniy*3}
  When $\eta$ satisfies $(\Hsw)$, the boundary $\yse(x)$ is the unique
  solution $y(x)$ of~(\ref{EDOy*_3})  on $\mathbb{R}^*_+$ that is
  increasing, such that $y(0^+)=0$ and $\min[y(x)-x,\delta(y(x)-x)+\gamma(r-\delta)y(x)]>0$ for all $x>0$.
\end{proposition}
\begin{proof}
   Let $y_1(x)$ and $y_2(x)$ denote two solutions of \eqref{EDOy*_3}
   satisfying the above hypotheses and $\tilde I(x)=\frac{y_1^{-1}(y_2(x))}{x}$.
  We have $\tilde I'(x)=\frac{1}{x}\left(\frac{F(x\tilde
      I(x),y_2(x))}{F(x,y_2(x))}-\tilde I(x)\right)$ where $$F(z,y)=2(y-z)[\delta(y-z) +\gamma(r-\delta)y
     ]+\gamma(1-\gamma)y^2\eta^2(y).$$Writing the estimations satisfied by $y_1$ (resp. $y_2$) at
$y_1^{-1}(y_2(x))$ (resp. $x$) one obtains $F(x\tilde I(x),y_2(x))>0$
(resp. $F(x,y_2(x))>0$). Moreover, since
$\partial_zF(z,y)=-2[2\delta(y-z)+\gamma(r-\delta)y]$ both $x\tilde I(x)$ and
$x$ belong to the interval
$(0,\frac{(2\delta+\gamma(r-\delta))y_2(x)}{2\delta})$ on which
$z \mapsto F(z,y_2(x))$ is decreasing. One easily concludes by the
same argument as in the proof of Proposition~\ref{uniy*}.
\end{proof}
\begin{proposition}\label{unix*3}
If $\sigma$ satisfies $(\Hsw)$ and $\max\left(r-\delta,\frac{(\delta-r)(\gamma\delta+(1-\gamma)r)}{(1-\gamma)\delta+\gamma r}\right)>\frac{(1-\gamma)\bar{\sigma}^2}{2}$,
    then $\xss(y)$ is the unique solution $x(y)$ of~(\ref{EDOx*_3}) on $\mathbb{R}^*_+$ that is
  increasing and such that
  $\exists \varepsilon>0$, $\forall y>0$, $\varepsilon y<x(y)<\min\left(1,\frac{(1-\gamma)\delta+\gamma r}{\delta}\right)y$.\end{proposition}
\begin{proof}
  By the convexity of $x \mapsto 1/x$, one has $r/(\gamma
  \delta+(1-\gamma)r)\le r(\gamma/\delta+(1-\gamma)/r)=(\gamma r+(1-\gamma)\delta)/\delta$.
  Therefore, using Remark~\ref{rem_1} for the first inequality, one deduces
  \begin{equation}
   \xss(y)<\frac{r}{\gamma
  \delta+(1-\gamma)r}y \le \frac{\gamma r+(1-\gamma)
  \delta}{\delta}y. \label{majxssss}
\end{equation}

Let $x(y)$ denote a solution of~\eqref{EDOx*_3} and $\tilde{I}(y)=\frac{x^{-1}(\xss(y))}{y}$, one has
$\tilde{I}'(y)=\frac{\tilde{I}(y)-1}{y}G(y)$ with $$G(y)=\frac{2[r(\tilde{I}(y)y^2-\xss(y)^2)+\gamma(r-\delta)\xss(y)^2]-\gamma(1-\gamma)\xss(y)^2\sigma^2(\xss(y))}{2[r(y-\xss(y))^2+\gamma(r-\delta)\xss(y)(y-\xss(y))]+\gamma(1-\gamma)\xss(y)^2\sigma^2(\xss(y))}.$$
 By Proposition \ref{compfront} and Example \ref{volconst},
 $\frac{a(\bar{\sigma})}{a(\bar{\sigma})-\gamma}y=x_{\bar{\sigma}}(y)\leq \xss(y)<y$, which implies that the denominator in the definition of $G$ is not greater than $\left(2\gamma^2[\frac{r}{a^2(\bar{\sigma})}-\frac{(r-\delta)^+}{a(\bar{\sigma})}]+\gamma(1-\gamma)\bar{\sigma}^2\right)\xss(y)^2$. If $r-\delta>\frac{(1-\gamma)\bar{\sigma}^2}{2}$ and $x(y)<y$, then $\xss(y)<x^{-1}(\xss(y))=y\tilde{I}(y)$ and 
$$\forall y>0,\;G(y)>\frac{\gamma a^2(\bar{\sigma})(2(r-\delta)-(1-\gamma)\bar{\sigma}^2)}{2\gamma^2[r-a(\bar{\sigma})(r-\delta)]+\gamma(1-\gamma)a^2(\bar{\sigma})\bar{\sigma}^2}>0.$$
If $\frac{(\delta-r)(\gamma\delta+(1-\gamma)r)}{(1-\gamma)\delta+\gamma
  r}>\frac{(1-\gamma)\bar{\sigma}^2}{2}$ and
$x(y)<\frac{(1-\gamma)\delta+\gamma r}{\delta}y$ then
$\frac{\delta}{(1-\gamma)\delta+\gamma
  r}\xss(y)<x^{-1}(\xss(y))=y\tilde{I}(y)$ and using the first inequality
in~\eqref{majxssss}, we get
$$\forall y>0,\;G(y)>\frac{\gamma a^2(\bar{\sigma})\left(\frac{2(\delta-r)(\gamma\delta+(1-\gamma)r)}{(1-\gamma)\delta+\gamma
  r}-(1-\gamma)\bar{\sigma}^2\right)}{2\gamma^2r+\gamma(1-\gamma)a^2(\bar{\sigma})\bar{\sigma}^2}>0.$$
In both cases, when $\tilde{I}(1)>1$ then $\forall y>0$,
$\tilde{I}(y)>1$ and for $y>1$, $\tilde{I}(y)-\tilde{I}(1)\geq
c(\tilde{I}(1)-1)\log(y)$ for some positive constant $c$. This
contradicts the inequality
$x^{-1}(\xss(y))<x^{-1}(y)<\frac{y}{\varepsilon}$ which holds as soon as
for all $y>0$, $x(y)>\varepsilon y$.
When $\tilde{I}(1)<1$, then for $y>1$, $\tilde{I}(y)-\tilde{I}(1)\leq
c(\tilde{I}(1)-1)\log(y)$ which contradicts the positivity of
$\tilde{I}$.
\end{proof}

\begin{theorem}\label{dualpart*3}
Let us assume  that
$\sigma$ and $\eta$ satisfy $(\Hsw)$. Then, the following assertions are equivalent:\\
(1) $\forall
  x,y>0, \ P_{\sigma}(x,y)=c_\eta(y,x)$.\\
(2) $\forall y>0$,
$\xss(y)^2\sigma^2(\xss(y))>\frac{2(1-\gamma)}{\gamma^2(2-\gamma)}[r(y-\xss(y))^2+\gamma(r-\delta)\xss(y)(y-\xss(y))]$
  and $\eta \equiv \tsigma$ where
  \begin{align*}\label{sigma_tilde}
    \tsigma^2(y)&=\frac{2[\delta(y-\xss(y))^2+\gamma(r-\delta)y(y-\xss(y))]}{\gamma y^2}\\&\times\frac{2[r(y-\xss(y))^2+\gamma(r-\delta)\xss(y)(y-\xss(y))]+\gamma(1-\gamma)\xss(y)^2\sigma^2(\xss(y))}{\gamma^2(2-\gamma)\xss(y)^2\sigma^2(\xss(y))-2(1-\gamma)[r(y-\xss(y))^2+\gamma(r-\delta)\xss(y)(y-\xss(y))]}.
  \end{align*}
If moreover $\max\left(r-\delta,\frac{(\delta-r)(\gamma\delta+(1-\gamma)r)}{(1-\gamma)\delta+\gamma
    r}\right)>\frac{(1-\gamma)\bar{\sigma}^2}{2}$, they are
also equivalent to\\
(3)
$\forall x>0,\yse(x)^2\eta^2(\yse(x))>\frac{2(1-\gamma)}{\gamma^2(2-\gamma)}[\delta(\yse(x)-x)^2+\gamma(r-\delta)\yse(x)(\yse(x)-x)]$,
$\yse(x)>\frac{\gamma\delta+(1-\gamma)r}{r}x$ and $\sigma(x)=\uteta(x)$ where
 \begin{align*}
    \uteta^2(x)&=\frac{2[r(\yse(x)-x)^2+\gamma(r-\delta)x(\yse(x)-x)]}{\gamma x^2}\\&\times\frac{2[\delta(\yse(x)-x)^2+\gamma(r-\delta)\yse(x)(\yse(x)-x)]+\gamma(1-\gamma)\yse(x)^2\eta^2(\yse(x))}{\gamma^2(2-\gamma)\yse(x)^2\eta^2(\yse(x))-2(1-\gamma)[\delta(\yse(x)-x)^2+\gamma(r-\delta)\yse(x)(\yse(x)-x)]}.
  \end{align*}\end{theorem}

Using~\eqref{majxssss}, the numerator in the first term of the r.h.s. of the equation
giving $\tsigma^2(y)$ is positive. Notice that in (3), the condition
$\yse(x)>\frac{\gamma\delta+(1-\gamma)r}{r}x$ that ensures the positivity of the
analogous term is satisfied as soon as
$\delta\leq r$.
\begin{proof}
   The proof is similar to the one of Theorem~\ref{dualpart}. For $(1)\Rightarrow
   (3)$, $\yse(x)>\frac{\gamma\delta+(1-\gamma)r}{r}x$ comes from $\yse \equiv
   {\xss}^{-1}$ and~\eqref{majxssss}.
For $(2)\Rightarrow (1)$, one remarks that according to \eqref{majxssss} and
Proposition~\ref{bound_boundaries},
$\forall x>0,\;{\xss}^{-1}(x)>\max(\frac{\delta x}{\gamma r+(1-\gamma)\delta},x)$ which, combined with Proposition~\ref{uniy*3}, ensures that
${\xss}^{-1}=y_{\tilde{\sigma}}$.\\
For $(3)\Rightarrow (1)$, since according to Proposition
\ref{compfront}, Example~\ref{volconst} and Remark \ref{rem_1},
$\frac{b(\bar{\eta})}{b(\bar{\eta})-\gamma}x\geq \yse(x)>\max\left(1,\frac{\delta
    }{\gamma r+(1-\gamma)\delta}\right)x$, one has
$\frac{(b(\bar{\eta})-\gamma)}{b(\bar{\eta})}y\leq {\yse}^{-1}(y)<\min\left(1,\frac{\gamma r+(1-\gamma)\delta}{\delta
    }\right)y$.\end{proof}
To obtain an analytical example of non-constant dual volatility
functions, we consider the same reciprocal boundaries as in~\cite{AJ} :
\[y^*(x)=x \frac{x+a}{bx+c}\;\mbox{ and }\;x^*(y)=\frac{1}{2}\left(
  by-a+\sqrt{(by-a)^2+4cy } \right)\;\mbox{ with }\;a,b,c>0.\]
Under some assumptions on $a,b,c$, these functions are the exercise
boundaries associated with explicit dual volatility functions.
\begin{proposition} Let us assume that
  $\max(c/a,b)<\min(1,\frac{r}{(1-\gamma)r+\gamma \delta})$,
  $\min(c/a,b)>1-\gamma$ and 
 $\frac{1-\gamma}{\gamma^2} [\max(\frac{1}{b},\frac{a}{c})-1] [
  r\max(\frac{1}{b},\frac{a}{c}) - ((1-\gamma)r+\gamma \delta)   ]<
  \max\left(r-\delta,\frac{(\delta-r)(\gamma\delta+(1-\gamma)r)}{(1-\gamma)\delta+\gamma r}\right)$.
  Then the volatility functions
 \begin{eqnarray*}
    \sigma(x) &=& \sqrt{ \frac{2}{\gamma}\frac{\big[x(1-b)+a-c \big] \big[x
        [r - b((1-\gamma)r+\gamma \delta) ]+a  [r - \frac{c}{a}((1-\gamma)r+\gamma \delta) ] \big]}{bx^2+2cx+ac +(\gamma-1)(bx+c)^2}} \\
    \eta(y)&=&\sqrt{\frac{2}{\gamma} \frac{[y-x^*(y)][\delta(y-x^*(y)) +\gamma(r-\delta)y][bx^*(y)^2+2cx^*(y)+ac]}{y^2[b(b+\gamma-1)x^*(y)^2+2c(b+\gamma-1)x^*(y)+ca(\frac{c}{a}+\gamma-1)]} }
\end{eqnarray*} 
are well defined and satisfy $(\Hs)$. Moreover we have $\yse\equiv y^*$ and
$\xss\equiv x^*$ and thus the duality holds: $\forall x,y>0, \
P_\sigma(x,y)=c_\eta(y,x)$.   
\end{proposition}
We have plotted in Figure~\ref{Dualgenfig} an example that illustrates
this duality result.
\begin{remark}
   Let us comment briefly the assumptions on the coefficients
   $a,b,c$. Under the
   second hypothesis, $\max(1/b,c/a)<(1-\gamma)^{-1}$ and therefore the third assumption will be automatically
satisfied if
\begin{eqnarray*}
&&  \frac{r}{\gamma^2} (\max(1/b,a/c)-1)  <r-\delta \text { when } r> \delta, \\
&&
\frac{r}{\gamma}(\max(1/b,a/c)  - ((1-\gamma) +\gamma \delta/r)   ) <\frac{(\delta-r)(\gamma\delta+(1-\gamma)r)}{(1-\gamma)\delta+\gamma r} \text { when } r< \delta.\end{eqnarray*}
When $r>\delta$, it is always possible to take $b$ and $a/c$ close
enough to $1$ so
that the equivalent condition $\max(1/b,a/c )< 1+ \gamma^2 (1-\delta/r)$ is satisfied. In the same manner,
if $\delta>r$, the first hypothesis writes $\min(1/b,a/c)>1-\gamma +\gamma
\delta/r$, and the third assumption will be always satisfied if one takes parameters
such that $1/b$ and $a/c$ are close enough to $1-\gamma +\gamma
\delta/r$. This is nonetheless compatible with the second assumption only if $1-\gamma<\frac{r}{r(1-\gamma)+\delta
\gamma}$, i.e. $\delta/r <\frac{2-\gamma}{1-\gamma}$. Otherwise,
there are no parameters $a,b,c$ that fulfill the three assumptions. Let us
remark incidentally that this condition is the same as the condition $\gamma-1 > -
\frac{r}{\delta-r}$ which appears in the Black-Scholes case (see Example~\ref{volconst}). 
\end{remark}

\begin{proof} {\it First step: let us check that the functions $\sigma$ and $\eta$ are well defined and
    satisfy $(\Hs)$.} The denominator in the definition of $\sigma$ is equal to
  $b(1+(\gamma-1)b)x^2+2c(1+(\gamma-1)b)x+ac(1+(\gamma-1)c/a)$: this is a second
  degree polynomial with positive coefficients because $\max(c/a,b)<1$. It is then  easy to check that $\sigma$ is well defined and satisfy $(\Hs)$
  using that $\max(c/a,b)<\min(1,\frac{r}{(1-\gamma)r+\gamma \delta})$. We get after some calculations
  \[ \eta(y^*(x))=\sqrt{\frac{2}{\gamma} \frac{[x(1-b)+a-c][ (
      \gamma r+(1-\gamma)\delta) - \delta b )x+(\gamma r+(1-\gamma)\delta - \delta
      \frac{c}{a})a
      ][bx^2+2cx+ac]}{(x+a)^2[b(b+\gamma-1)x^2+2c(b+\gamma-1)x+ca(\frac{c}{a}+\gamma-1)]} }\]
{From} the first hypothesis and the argument given at the beginning of the proof of Proposition~\ref{unix*3},
one obtains $\max(c/a,b)<\min(1,\frac{\gamma r+(1-\gamma)\delta)}{\delta})$.
Using the second hypothesis and the one-to-one onto
property of the function $y^*$, we deduce
that $\eta$ is also well
defined and satisfies $(\Hs)$ .

{\it Second step:}
We easily check that $x^*$ and $y^*$ respectively solve the
ODEs~(\ref{EDOx*_3}) and~(\ref{EDOy*_3}). Since
$\max(c/a,b)<\min(1,\frac{\gamma r+(1-\gamma)\delta)}{\delta})$, we have  $y^*(x)>x$ and
$\delta(y^*(x)-x)+\gamma(r-\delta)y^*(x)>0$. Proposition~\ref{uniy*3}
then ensures that
$y^*\equiv\yse$. Since $bx^2+2cx+ac>(bx+c)^2$,
we have
\begin{align*}
   \sigma^2(x)&\le \frac{2}{\gamma^2} \frac{x(1-b)+a-c}{bx+c} \frac{x
        [r - b((1-\gamma)r+\gamma \delta) ]+a  [r -
        \frac{c}{a}((1-\gamma)r+\gamma \delta) ]}{bx+c}\\
&\leq \frac{2}{\gamma^2}\left[\max(\frac{1}{b},\frac{a}{c})-1\right] \left[
  r\max(\frac{1}{b},\frac{a}{c}) - ((1-\gamma)r+\gamma \delta)\right ]  \\
&\leq \frac{2}{1-\gamma}\max\left(r-\delta,\frac{(\delta-r)(\gamma\delta+(1-\gamma)r)}{(1-\gamma)\delta+\gamma r}\right).
\end{align*}
By Proposition~\ref{unix*3}, we
conclude that $x^*\equiv\xss$. \end{proof}

\subsection{The theoretical calibration procedure}

In \cite{AJ}, the calibration issue was the practical motivation for our interest in
Call-Put duality. Even if, in the present framework, calibration is
purely theoretical since the payoff $\phi(x,y)$ is not traded, we are
going to explain shortly how to recover the local volatility function
from the perpetual prices of options. More precisely, let us suppose
that we observe for all $K>0$ the market price $p(K)$ of the
American security with payoff $\phi(.,K)$, with either $\phi(x,y)=((y-x)^+)^\gamma$
and $\gamma \in (0,1]$ or $\phi(x,y)=(\psi_y(y)-\psi_x(x))^+$, with $\psi_y$ and
$\psi_x$ satisfying the assumptions mentioned before. We also assume that either $\max\left(r-\delta,\frac{(\delta-r)(\gamma\delta+(1-\gamma)r)}{(1-\gamma)\delta+\gamma
    r}\right)>\frac{(1-\gamma)\bar{\sigma}^2}{2}$ or
$\psi_x$ satisfy~(\ref{Hyp_psix}) so that we have equivalence between the three
conditions in Theorems ~\ref{dualpart*3} or~\ref{dualpart}. We denote by 
$x_0$ the current value of the stock, and we suppose that there is a
function $\sigma$  satisfying $(\Hsw)$ such that $\forall K>0$,
$p(K)=P_{\sigma}(x_0,K)$ and that
$\tsigma$ is well defined and satisfy $(\Hsw)$. Within this framework, as in the
call-put case, we are able to get $(\sigma(x),0< x\le x_0)$.

Indeed, let us define $Y=\inf\{K>0, p(K)= \phi(x_0,K) \}$. Since $p(K)=c_{\tsigma}(K,x_0)$, we have 
$Y=y^*_{\tsigma}(x_0)$ and \[\forall K<Y, \ \frac{K^2 \tsigma(K)^2}{2}p''(K)+K(\delta-r)p'(K)-\delta p(K)=0 \]
We deduce then $\forall K \le Y,  \tsigma(K)=\frac{1}{K} \sqrt{ \frac{2(\delta
    p(K)+K(r-\delta)p'(K))}{p''(K)}}$ because
$p''(K)=\partial^2_Kc_{\tsigma}(K,x_0)=\frac{\phi(x_0,Y)}{g(Y)}g''(K)>0$
for $K<Y$ using~(\ref{convex_fg}). Then, we get $(y^*_{\tsigma}(x),0< x\le x_0)$ solving backward either
(\ref{EDOy*_3}) or (\ref{EDOy*_2}) starting from $y^*_{\tsigma}(x_0)=Y$. Finally, we get $(\sigma(x),0< x\le x_0)$ thanks
to Theorem~\ref{dualpart*3} or~\ref{dualpart}, using that $\sigma(x)=\tilde{\utsigma}(x)$.

Now, to formalize our
calibration result, we introduce the set $$\Sigma=\{\sigma\text{ satisfying } (\Hsw)\mbox{ s.t. }\tsigma\mbox{ is
  well defined and satisfies } (\Hsw)\}.$$
\begin{proposition}
  Under the above assumptions on $\phi$, $r$ and $\delta$, for $\sigma_1,\sigma_2\in\Sigma$,
$$\forall K>0,\;
P_{\sigma_1}(x_0,K)=P_{\sigma_2}(x_0,K)\;\Leftrightarrow\;{\sigma_1}\big|_{(0,x_0]}
\equiv {\sigma_2}\big|_{(0,x_0]}\mbox{ and }y^*_{\tilde{\sigma}_1}(x_0)=y^*_{\tilde{\sigma}_2}(x_0).$$
\end{proposition}

\begin{proof}
The necessary condition is a consequence of the above calibration
procedure. To check the sufficient condition, we consider $\sigma_1$ and $\sigma_2$ in $\Sigma$ such that
  \[ \forall x \le x_0, \sigma_1(x)=\sigma_2(x)\text{
    and } y^*_{\tsigma_1}(x_0)=y^*_{\tsigma_2}(x_0)=Y.\]
  On the one hand, we have $x^*_{\sigma_1}(Y)=x^*_{\sigma_2}(Y)=x_0$, and thus
  $x^*_{\sigma_1}(y)=x^*_{\sigma_2}(y)$ for $y \le Y$ since they solve the same ODE.
  Therefore, using either Theorem~\ref{dualpart*3} or Theorem~\ref{dualpart}, one gets
  $\tsigma_1(y)=\tsigma_2(y)$ for $y\le Y$. On
  the other hand, the smooth fit principle gives
  $\frac{g'_{\tsigma_1}(Y)}{g_{\tsigma_1}(Y)}=\frac{\partial_y\phi(x_0,Y)}{\phi(x_0,Y)}=\frac{g'_{\tsigma_2}(Y)}{g_{\tsigma_2}(Y)}$.
  The set of solutions to
  $\frac{1}{2}y^2\tsigma_1^2(y)g''(y)+(\delta-r)yg'(y)-\delta g(y)=0$ on $(0,Y]$ is a
  two-dimensional vectorial space, and by the previous equality,
  $g_{\tsigma_1}$ and $g_{\tsigma_2}$ are proportional on $(0,Y]$. Therefore,
  we have for $0<K\le Y$,
  $P_{\sigma_1}(x_0,K)=\phi(x_0,Y)\frac{g_{\tsigma_1}(K)}{g_{\tsigma_1}(Y)}=\phi(x_0,Y)\frac{g_{\tsigma_2}(K)}{g_{\tsigma_2}(Y)}
  =P_{\sigma_2}(x_0,K)$, and $P_{\sigma_1}(x_0,K)=\phi(x_0,K)=P_{\sigma_2}(x_0,K)$ for
  $K\ge Y$.
\end{proof}
Like in  Proposition~5.1 \cite{AJ}, we can get an analogous calibration of the
complementary upper part of the local volatility function to the
perpetual prices of the ``Call'' options with payoff $\phi(K,x_0)$ by exchanging the roles of $\eta$ and $\sigma$, and of $r$
and $\delta$.

\end{document}